\documentclass[12pt,a4paper]{article}
\usepackage{latexsym,amsfonts,amssymb,epsfig,euscript}
\usepackage{roundop,matstyle}
\usepackage{color}
\usepackage{url}
\definecolor{red,green,blue}{rgb}{0.5,0.5,0.7}
\definecolor{cyan,magenta,yellow,black}{cmyk}{0.5,0.5,0.5,1}

\def\R2n{{\bf {R}}^{2n}}

\def\R{{\bf {R}}}

\def\ux{\underline{x}}
\def\ox{\overline{x}}
\def\ud{\underline{d}}
\def\od{\overline{d}}

\def\ty{\tilde{y}}

\def\tx{\tilde{x}}

\def\tX{\tilde{X}}

\def\tD{\tilde{D}}
\def\tQ{\tilde{Q}}
\def\tLambda{\tilde{\Lambda}}

\def\hf{\hat{f}}

\def\oX{\overline{X}}
\def\of{\overline{f}}
\def\oy{\overline{y}}
\def\on{\overline{n}}
\def\ore{\overline{r}}
\def\uy{\underline{y}}

\def\ux{\underline{x}}

\def\tlambda{\tilde{\lambda}}

\def\uf{\underline{f}}
\def\trace{\mbox{\rm trace}}
\def\Diag{\mbox{\rm Diag}}
\def\diag{\mbox{\rm diag}}
\def\rank{\mbox{\rm rank}}
\def\svec{\mbox{\rm svec}}
\def\EX{\EuScript{X}}
\def\EK{\EuScript{K}}
\def\EY{\EuScript{Y}}
\def\EM{\EuScript{M}}

\def\ohf{\overline{f}}
\def\uD{\underline{D}}

\newtheorem{theorem}{Theorem}[section]
\newtheorem{coro}{Corollary}[section]
\newtheorem{lemma}{Lemma}[section]
\newtheorem{prop}{Proposition}[section]

\title{Guaranteed Accuracy for Conic Programming Problems in Vector Lattices}
\author{Christian ~Jansson}

\begin{document}

\maketitle

\vspace*{0.5cm} Institute for Reliable Computing, Technical
University Hamburg--Harburg, Schwarzenbergstra{\ss}e 95, 21071
Hamburg, Germany, e-mail: jansson@tu-harburg.de, Fax: ++49 40
428782489.

\vspace*{0.5cm} \noindent{\bf{Keywords}}: Linear programming,
semidefinite programming, conic programming, convex programming,
combinatorial optimization, rounding errors, ill-posed problems,
interval arithmetic, branch-bound-and-cut
\\
\vspace*{0.5cm}

\begin{abstract}
This paper presents rigorous forward error bounds for linear conic
optimization problems. The error bounds are formulated in a quite
general framework; the underlying vector spaces are not required
to be finite-dimensional, and the convex cones defining the
partial ordering are not required to be polyhedral. In the case of
linear programming, second order cone programming, and
semidefinite programming specialized formulas are deduced yielding
guaranteed accuracy. All computed bounds are completely rigorous
because all rounding errors due to floating point arithmetic are
taken into account. Numerical results, applications and software
for linear and semidefinite programming problems are described.
\end{abstract}

\section{Introduction}

In this paper forward error bounds for the optimal value of linear
conic optimization problems as well as certificates of feasibility
and infeasibility are presented, including the discussion of
rounding-off errors and details of implementation. These rigorous
bounds aim to prove how accurate the approximate results computed
by any conic solver are. The underlying vector spaces are in
general  infinite-dimensional, that is the bounds are developed in
the framework of functional analysis.

Forward and backward error analysis together with a detailed
discussion of rounding-off errors and condition numbers for matrix
problems were first described in the outstanding papers published
sixty years ago by von Neumann and Goldstine \cite{NeuGo47}  and
Turing \cite{Tur48}. Turing writes:
\begin{quote}
Error estimates can be of two kinds. We may wish to know how
accurate a certain result is, and be willing to do some additional
computation to find out. A different kind of estimate is required
if we are planning calculations and wish to know whether a given
method will lead to accurate results. In the former case we do not
care what quantities the error is expressed in terms of, provided
they are reasonably easily computed. With these estimates we wish
to be absolutely sure that the error is within the range stated,
but at the same time not to state a range which is very much
larger than necessary. With the second type of estimate, the error
is preferably expressed in terms of quantities whose meaning is
sufficiently familiar that the general run of values involved may
at least be guessed at.
\end{quote}
Particularly, forward error bounds for the inverse of a matrix
including a discussion of the effects of rounding-off errors are
presented there. Today one would speak in this context of verified
or rigorous error bounds, and thus these two papers can be viewed
as the pioneering work in the field of verification methods, a
part of numerical analysis. Forward error bounds are propagated in
interval arithmetic; see the textbooks Alefeld and Herzberger
\cite{AlHe83}, Moore \cite{Mo79}, and Neumaier \cite{Ne90},
\cite{Neu2001}. But also in other areas the interest in rigorous
forward error bounds is growing. Parlett \cite{Par01}, for
example, remarks in relation to the numerical accuracy of
eigenvalue problems:
\begin{quote}
For some of us, however, it has taken nearly 40 years to realize
that backward stability is not enough.
\end{quote}
Also Trefethen writes in \cite{Tre06} about the future of
Numerical Analysis
\begin{quote}
I expect that most of the numerical computer programs of 2050 will
be 99\% intelligent and just 1\% actual ``algorithm'' if such a
distinction makes sense. Hardly anyone will know how they work,
but they will be extraordinarily powerful and reliable, and will
often deliver results of guaranteed accuracy.
\end{quote}

Linear conic optimization refers to problems with a linear
objective function and linear constraints where the variables are
restricted to a cone. In general these problems are non-smooth.
Linear programming, quadratically convex programming, second order
cone programming and semidefinite programming are special cases.
Since each convex problem can be described equivalently as a
linear conic problem, the latter provides a universal form of
convex programming (see Nemirovski \cite{Nem96},
\cite{NestNem94}). Thus not surprisingly, a large variety of
applications of conic programming are known from areas like system
and control theory, combinatorial optimization, signal processing
and communications, machine learning, quantum chemistry, and many
others. For an elaborate bibliography the reader is referred to
Wolkowicz \cite{Wolko}.

Nesterov and Nemirovski \cite{NestNem94} have shown that
self-concordant barriers apply to many conic problems yielding
polynomial time interior point methods. Renegar \cite{Ren94} has
investigated the sensitivity of infinite-dimensional conic
optimization problems, and in \cite{Renegar95} he analyzed
interior point methods. He introduced a condition number for conic
optimization, which is a generalization of the condition numbers
defined by von Neumann, Goldstine, and Turing. This
\emph{condition number} is the scale-invariant reciprocal of the
smallest  data perturbation that will render the perturbed data
instance either primal or dual infeasible. It is used in
sensitivity analysis and moreover can be viewed as a problem
instance size of conic optimization problems yielding important
results in complexity theory. A problem is called \emph{ill-posed}
if this condition number is infinite, that is the distance to
primal or dual infeasibility is zero.

One of Renegar's main results is that the sensitivity of the
optimal solutions and the optimal value can be bounded by the
condition, and especially he proved that the bounds for the
optimal value depend cubically on the inverses of the relative
distances to primal and dual infeasibility. Renegar shows that
this bound cannot be improved in general. For an ill-posed problem
this result means that there exist arbitrarily small perturbed
data instances such that the difference between the optimal value
of the original problem and the perturbed problem is arbitrarily
large, but the optimality conditions for the perturbed problem
almost coincide with the optimality conditions for the original
problem. Since conic solvers are terminated when the optimality
conditions are satisfied approximately it cannot be distinguished
between the optimal values of the original and the perturbed
problem in the case of ill-conditioned or ill-posed problems. A
consequence is that the noise introduced by floating point
arithmetic may occasionally yield to wrong termination and
nonsensical computational results.

In this paper we show how certain weak boundedness qualifications
on $\varepsilon$-optimal solutions can be used to compute rigorous
forward error bounds for the exact optimal value, also for
ill-conditioned or even for ill-posed problems. Such
qualifications and even more restrictive assumptions, like certain
smoothness properties, are customary when solving ill-posed
problems with regularization methods. It need not to be assumed
that Slater's constraint qualifications are fulfilled. The
rigorous error bounds provide more safety for conic optimization
problems, and they provide rigorous results in branch-and-bound
algorithms for global and combinatorial optimization problems.
Another application are computer-assisted proofs where it is
mandatory to control all rounding-off errors (see for example
Neumaier \cite{Ne06} and Rump \cite{Ru05a}). It should be made
clear that we do not investigate regularization methods. In this
paper we assume that approximations computed by some conic solver
(with or without regularization)  are given, and these
approximations are then used for computing the error bounds. It is
of particular importance that the computation of the error bounds
can be done outside the code of any imaginable solver as a
reliable postprocessing routine, providing a correct output for
the given input. Especially, we show for some combinatorial
problems how branch-and-bound algorithms can be made safe, even if
ill-posed relaxations are used. Numerical results for some
ill-posed and ill-conditioned problems are included.

Ill-conditioned and ill-posed problems are not rare in practice,
they occur even in linear programming. Ord{\'{o}}{\~{n}}ez and
Freund \cite{OrdFre03} stated that 71\% of the lp-instances in the
NETLIB Linear Programming Library \cite{Netlib} are ill-posed.
This library contains many industrial problems. Recently Freund,
Ord{\'{o}}{\~{n}}ez and Toh 2006 \cite{FreOrdToh06} have shown
that $32$ out of $85$ problems of the SDPLIB are ill-posed.

The presented results in this paper formalize a viewpoint which
apparently has not been made in conic programming. They can be
viewed as an extension of methods for linear programming (
\cite{Ja04a} and Neumaier and Shcherbina \cite{Neumaier5}), and
for smooth convex programming (see \cite{Ja04}) to ill-conditioned
and ill-posed non-smooth problems using the framework of
functional analysis.

The paper is organized as follows. After introducing some notation
and basic definitions in Section 2, we consider in the next
section conic optimization problems. Then in Section 4 verified
lower and upper bounds of the optimal value in the
infinite-dimensional case are presented, and applied to
finite-dimensional linear programming problems. Sections 5 and 6
are devoted to error bounds for second order cone and semidefinite
programming, respectively. Then in Section 7 we investigate conic
optimization problems with block structured variables. In Section
8 verified certificates of infeasibility are presented, and in the
following section we will focus on some applications in
combinatorial optimization. Section 10 contains numerical results
for the NETLIB Linear Programming Library (obtained by the C++
software package Lurupa \cite{Keil05}) and for the SDPLIB
benchmark problems (obtained by the MATLAB software package VSDP
\cite{Ja06a}). Finally, some conclusions are given.

\section{Notation and Preliminaries}

Let $\EX$ be a real vector space equipped with a norm $\|.\|$, and
let $\EK \subseteq \EX$ be a \emph{convex cone}, i.e.
\begin{equation} \label{eq2.1}
\EK + \EK \subseteq \EK, \; \alpha \EK \subseteq \EK \; \mbox{for}
 \; \alpha \in \R_+,
\end{equation}
where $\R_+$ denotes the set of nonnegative real numbers. A convex
cone $\EK$ induces a \emph{partial ordering}
\begin{equation} \label{eq2.2}
x \le y \quad \Leftrightarrow \quad y-x \in \EK,
\end{equation}
which is a transitive and reflexive binary relation on $\EX$
compatible with addition and scalar multiplication:
\begin{equation} \label{eq2.3}
x \le y, \; u \le v, \; \alpha \in \R_+ \quad \Rightarrow \quad
x+u \le y+v \; \mbox{and} \; \alpha x \le \alpha y.
\end{equation}
Conversely, each partial ordering determines a convex cone, namely
the \emph{positive cone}
\begin{equation} \label{eq2.4}
\EK := \{x \in \EX: \; x \ge 0\}.
\end{equation}
A vector space $\EX$ equipped with a partial ordering is called a
\emph{partially ordered vector space}. A partial ordering is
called \emph{antisymmetric}, if
\[
x \le y, \; y \le x \quad \Rightarrow \quad x = y.
\]
It can be proved that antisymmetric partial orderings correspond
to \emph{pointed cones}, i.e.
\[
\EK \cap(-\EK) = \{0\}.
\]

If not explicitely mentioned we do not assume that the partial
ordering is antisymmetric. Given a partial ordering the set
\begin{equation} \label{eq2.5}
[\ux,\ox] := \{x \in \EX: \; \ux \le x \le \ox\} = (\ux+\EK) \cap
(\ox - \EK)
\end{equation}
is called an \emph{interval}. For a subset $\EM$ of a partially
ordered vector space $\EX$ a vector $\ux$ is called a \emph{lower
bound} of $\EM$, if $\ux \le m$ for all $m \in \EM$, and in this
case we write $\ux \le \EM$. The lower bound $\ux$ is called
\emph{infimum} of $\EM$ if every other lower bound $\uy$ of $\EM$
satisfies $\uy \le \ux$. Analogously, \emph{upper bounds} and
\emph{supremum} are defined. $\EX$ is said to be a \emph{vector
lattice} for the partial ordering $\le$ if  for all $x,y \in \EX$
the supremum $\sup\{x,y\}$ and the infimum $\inf\{x,y\}$ exists
and is contained in $\EX$, respectively. In a vector lattice the
operations $x^+:= \sup\{x,0\}$, $x^-:=\inf\{x,0\}$ and
$|x|:=\sup\{x,-x\}$ are defined, and the properties $|x|=x^+-x^-,
\; x=x^++x^-, \;$ $|x|=0$ iff $x =0$, $|\lambda x| =
|\lambda|\;|x|$ for real $\lambda$, and $|x+y| \le |x|+|y|$ are
satisfied.

Let $\EX^*$ denote the \emph{dual space} of $\EX$, that is the
space of continuous linear functionals endowed with the operator
norm. The set $\EK^*$ of all positive linear functionals, i.e.
\begin{equation} \label{eq2.6}
\EK^* = \{y \in \EX^*: \; \langle y,x \rangle := y(x) \ge 0 \;
\mbox{for all } x \in \EK\},
\end{equation}
is a convex cone in $\EX^*$ defining a partial ordering in the
dual space.

The basic properties and relations for vector lattices as well as
examples can be found in Birkhoff \cite{Bir48} and Bourbaki
\cite{Bou52}; see also Peressini \cite{Pere67} and Schaefer
\cite{Schae74}. We use the same notation $\|.\|$ and $\le$ for all
norms and partial orderings. It will always be clear from the
context which norm and which cone is referred to. Hence, if $x \in
\EX$ then $x \ge 0$ means $x \in \EK$, and if $y \in \EX^*$ then
$y \ge 0$ denotes $y \in \EK^*$. Observe that we do not write
$y^*$ for a continuous linear functional in $\EX^*$ because from
the position in $\langle y,x \rangle$ the meaning is clear, and we
can omit the star. This notion is closely related to Hilbert
spaces and the Theorem of Riesz which states that the continuous
linear functions can be represented by the inner product $\langle
y,x \rangle$ where $y$ is a vector in the Hilbert space.

In the following some illustrative and well-known examples of
normed vector lattices are shown. The real finite dimensional
space $\EX = \R^n$ equipped with the Euclidean inner product and
the Euclidean norm $\|.\|$ can be ordered by the convex cone
\begin{equation} \label{eq2.7}
\EK:= \R^{n}_+=\{x \in \R^n: \; x_i \ge 0 \quad \mbox{for } i=1,
\ldots, n\}.
\end{equation}
This cone is \emph{self-dual} (i.e. $\EK = \EK^*$) and implies the
lattice operations
\begin{equation} \label{eq2.8}
x_i^+ = \max\{0, x_i\}, \; x_i^- = \min\{0,x_i\}, \;
|x_i|=x_i^+-x_i^-
\end{equation}
for $i = 1, \ldots, n$. This vector lattice is used in
\emph{linear programming} (LP).

In \emph{second order cone programming} (SOCP) the same normed
space $\EX=\R^n$ is equipped with the partial ordering defined by
the convex \emph{ice-cream} or \emph{Lorenz cone}
\begin{equation} \label{eq2.9}
\EK:= \left\{x={x_{:} \choose x_n} \in \R^n: \; x_n \ge \|x_{:}\|
\right\},
\end{equation}
where $x_{:}:= (x_1, \ldots, x_{n-1})^T$. This cone is also
self-dual and further properties are described in Section 5.

In \emph{semidefinite programming} (SDP) the real linear space
$\EX$ is $\R^{n(n+1)/2}$, which is identified with the set of real
symmetric $n \times n$ matrices $X$. The inner product of $X,Y$ is
defined by $\langle X,Y \rangle := \trace X^TY =
\Sigma_{ij}X_{ij}Y_{ij}$, and the induced norm $\|X\|:=(\trace X^T
 X)^{\frac{1}{2}}$ is the \emph{Frobenius norm}.

The space $\EX=\R^{n(n+1)/2}$ is a Hilbert space, thus self-dual,
and it is equipped with the self-dual cone of positive
semidefinite matrices
\begin{equation} \label{eq2.11}
\EK:= S_+^n= \{X \in \EX: \; X \mbox{ is positive semidefinite}\}.
\end{equation}
Using the eigenvalue decomposition $X=Q^T \Lambda Q$ of a real
symmetric matrix it follows that
\begin{equation} \label{eq2.12}
X^- = Q^T \Lambda^-Q, \quad X^+ = Q^T \Lambda^+Q, \quad |X|=Q^T
|\Lambda|Q,
\end{equation}
where $\Lambda^-$, $\Lambda^+$, and $|\Lambda|$ denote the
diagonal matrices with nonpositive, nonnegative, and modulus of
the eigenvalues of $X$ on the diagonal, respectively.

Occasionally, it is useful to represent symmetric matrices $X$ as
column vectors $x$ by using the \emph{svec} operator:
\begin{equation} \label{eqz1}
x:= \svec(X) := (X_{11}; \sqrt{2} X_{21}; \ldots ; \sqrt{2}
X_{n1}; X_{22}; \sqrt{2} X_{32}; \ldots ; X_{nn}).
\end{equation}
Here we follow the convention of MATLAB and use ``;'' for
adjoining vectors in a column. Then the inner product between
symmetric matrices $X$ and $Y$ is the usual inner product, that is
\begin{equation} \label{eqz2}
\langle X,Y \rangle  = x^Ty.
\end{equation}

We also use the notation $x \in \EK$ and $x \le y$ if the
corresponding symmetric matrices $X$ and $Y$ such that $x =
\svec(X)$ and $y = \svec(Y)$ have these properties.

For any compact Hausdorff space $\Omega$, the vector space
$\EX:=C(\Omega)$ of real-valued functions is a normed vector
lattice with norm
\[
\|f\|_{C(\Omega)}:= \sup\limits_{x \in \Omega} \{|f(x)|\}
\]
and ordering cone
\[
\EK:= \{f \in C(\Omega): \; f(x) \ge 0 \quad \mbox{for all } x \in
\Omega\}.
\]
Finally, we mention $L_p(\Omega)$, the vector space of
Lebesgue-integrable functions $f:\Omega \rightarrow \R$, where
$\Omega \subseteq \R^n$, $1 \le p < \infty$. This space is
equipped with the norm
\[
\|f\|_p := (\int\limits_{\Omega} |f(x)|^p dx)^{\frac{1}{p}},
\]
and can be partially ordered by the cone
\[
\EK:=\{f \in L_p(\Omega): \; f(x) \ge 0 \quad \mbox{almost
everywhere on } \Omega\}.
\]
This yields a normed vector lattice, which is of interest in the
case of Volterra and Fredholm type equations.

\section{Conic Optimization Problems}

We study rigorous error bounds for the \emph{conic optimization
problem in standard form}
\begin{equation} \label{eq3.1}
\mbox{minimize } \langle c,x \rangle \quad \mbox{s.t.} \quad Ax=b,
\; x \in \EK,
\end{equation}
where $\EX$ is a real normed vector space, $\EK \subseteq \EX$ is
a convex cone, $c \in \EX^*$, $\EY$ is a real normed vector space,
$A$ denotes a continuous linear operator from $\EX$ to $\EY$, and
$b \in \EY$. With $\hf_p$ we denote the primal optimal value,
where $\hf_p:= + \infty$ if the problem is infeasible. Many
interesting examples of optimization problems can be formulated in
this framework. In the following some familiar facts are
described.

The \emph{Lagrangian function} of problem (\ref{eq3.1}) has the
form
\begin{equation} \label{eq3.2}
L(x,y):= \langle c,x \rangle + \langle y,b-Ax \rangle,
\end{equation}
where $y \in \EY^*$. The optimization problem
\begin{equation} \label{eq3.3}
\inf\limits_{x \in \EK} \; \sup\limits_{y \in \EY^*} L(x,y)
\end{equation}
is equivalent to (\ref{eq3.1}). Indeed, if $b-Ax = 0$ then
$\langle y,b-Ax \rangle =0$ for each $y \in \EY^*$, and the
supremum of $L(x,y)$ is equal to $\langle c,x \rangle$. In the
case where $b-Ax \neq 0$ there is some $y$ with $\langle y,b-Ax
\rangle
> 0$, and hence the supremum is infinite.

Obviously, the Lagrangian satisfies $L(x,y)=\langle y,b \rangle +
\langle -A^*y+c,x\rangle$ where $A^*$ is the \emph{adjoint
operator}. By exchanging in (\ref{eq3.3}) infimum and supremum we
obtain the dual problem
\begin{equation} \label{eq3.4}
\sup\limits_{y \in \EY^*} \inf\limits_{x \in \EK} L(x,y)
\end{equation}
with optimal value $\hf_d$. Since exchanging $\inf$ and $\sup$
always produces a lower bound, \emph{weak duality} holds, that is
$\hf_d \le \hf_p$. Because $\inf\limits_{x \in \EK}L(x,y) = -
\infty$ whenever $-A^*y+c \notin \EK^*$ the dual problem can be
written equivalently in the form
\begin{equation} \label{eq3.5}
\mbox{maximize } \langle y,b \rangle \; \mbox{s.t.} \; -A^*y+c \in
\EK^*, \; y \in \EY^*.
\end{equation}
We set $\hf_d := - \infty$, if the dual problem is infeasible.

Let $x$ be primal feasible, and let $y$ be dual feasible, then
\begin{equation} \label{eq3.6}
\langle c,x \rangle = \langle c,x \rangle + \langle y,b-Ax \rangle
= \langle -A^* y+c, x \rangle + \langle y,b \rangle \ge \langle
y,b \rangle,
\end{equation}
and hence equality holds iff the \emph{complementarity condition}
\begin{equation} \label{eq3.7}
\langle -A^*y+c,x \rangle =0
\end{equation}
are fulfilled. This condition means that the feasible pair $x,y$
is a saddle point of the Lagrangian. Moreover, it follows that
there is no duality gap between the primal and the dual problem,
and both problems have optimal solutions if and only if there
exists a primal and a dual feasible solution fulfilling  the
complementarity conditions. In other cases where such  primal dual
feasible pairs do not exist strong duality may be not fulfilled.

Duality theory is central to the study of optimization. First,
algorithms are frequently based on duality (like primal-dual
interior point methods), secondly they enable to check whether a
given feasible point is optimal, and thirdly it allows to perform
a sensitivity analysis. For more results on duality theory in the
infinite-dimensional case see  for example Renegar \cite{Ren94},
Rockafellar \cite{Rock70}, and the literature cited there.

\section{Lower and Upper Bounds for the Optimal Value}

This section is elementary but important for understanding both,
the basic ideas behind rigorous forward error bounds and
implementations. It turns out that for computing these error
bounds only approximate primal and dual solutions $\tx$ and $\ty$
are required. Further assumptions about the accuracy of the
approximations are not necessary; they need to be neither primal
nor dual feasible. If the accuracy is poor, however, then the
error bounds cause overestimation.

The cones $\EK$ and $\EK^*$ create partial orderings for the
vector spaces $\EX$ and $\EX^*$, respectively. {\it We assume in
this paper that for subsets of these partially ordered vector
spaces there exist lower and upper bounds. If the existence of the
infimum or the supremum is necessary we mention this explicitely}.
Note that all vector lattices satisfy this assumption.

We start with a simple result concerning bounds for linear
functionals.

\begin{lemma} \label{lem4.1}
Assume that $x$, $\ox \in \EK$, $x \le \ox$, and let  $d, \ud^-
\in \EX^*$ with $\ud^- \le  \{d,0\}$. Then
\begin{equation} \label{eq4.1}
\langle d,x \rangle \ge \langle \ud^-, \ox \rangle \quad
\mbox{and} \quad \langle \ud^-, \ox \rangle \le 0.
\end{equation}
\end{lemma}

{\sc Proof.} Since $0 \le d-\ud^-$ and $x \ge 0$ it is
\[
0 \le \langle d-\ud^-, x \rangle = \langle d,x \rangle - \langle
\ud^-,x \rangle.
\]
Hence $\langle d,x \rangle \ge \langle \ud^-,x\rangle$. Since
$-\ud^- \ge 0$ and $\ox -x \ge 0$ the linearity of the functional
$-\ud^-$ yields
\[
0 \le \langle -\ud^-, \ox -x \rangle = \langle \ud^-, x \rangle -
\langle \ud^-, \ox \rangle,
\]
which immediately implies (\ref{eq4.1}). \hfill $\square$

The dual version of this lemma is:
\begin{lemma} \label{lem4.2}
Assume that $x, \ux^- \in \EX$ with $\ux^- \le \{x,0\}$, and let
$d, \od \in \EK^*$ with $d \le \od$. Then
\begin{equation} \label{eq22}
\langle d,x \rangle \ge \langle \od, \ux^- \rangle \quad
\mbox{and} \quad \langle \od,\ux^- \rangle \le 0.
\end{equation}
\end{lemma}

{\sc Proof.} Since $0 \le x-\ux^-$ and $d \ge 0$ it is
\[
0 \le \langle d,x-\ux^- \rangle = \langle d,x \rangle - \langle
d,\ux^- \rangle.
\]
Hence $\langle d,x \rangle \ge \langle d, \ux^- \rangle$. Since
$-\ux^- \ge 0$ and $\od-d >0$ it follows that
\[
0 \le \langle \od-d, -\ux^- \rangle = \langle d,\ux^-\rangle -
\langle \od, \ux^- \rangle
\]
which implies (\ref{eq22}). \hfill $\square$

\begin{lemma} \label{lem4.3}
Assume that $\EX$ and $\EX^*$ are normed vector lattices. Assume
that for $x, \ox \in \EX$ and $d, \od \in \EX^*$ the inequalities
\begin{equation} \label{eq23}
|x| \le \ox \quad \mbox{and} \quad |d| \le \od
\end{equation}
are satisfied. Then
\[
|\langle d,x \rangle| \le \langle \od, \ox \rangle \le \|\od\| \;
\|\ox\|.
\]
\end{lemma}
{\sc Proof.} Since $|x| = \sup\{x, -x\} \le \ox$ it follows that
\[
-\ox \le x \le \ox.
\]
Analogously we obtain
\[
-\od \le d \le \od.
\]
The inequalities $\ox - x \ge 0$ and $\od - d \ge 0$ imply
\[
0 \le \langle \od-d, \ox-x \rangle = \langle \od, \ox \rangle -
\langle \od, x \rangle - \langle d, \ox \rangle + \langle d,x
\rangle,
\]
and the inequalities $x + \ox \ge 0$ and $d + \od \ge 0$ imply
\[0
\le \langle d+\od, x+\ox \rangle = \langle d,x \rangle + \langle
d,\ox \rangle + \langle \od, x \rangle + \langle \od, \ox \rangle.
\]
Adding both inequalities yields
\[
0 \le 2(\langle \od, \ox \rangle + \langle d,x \rangle ),
\]
and therefore
\[
-\langle d,x \rangle \le \langle \od, \ox \rangle.
\]
Because of the symmetry of $x$ and $-x$ in the definition
$|x|=\sup\{x, -x\}$ it follows that
\[
\langle d,x \rangle = -\langle d,-x \rangle \le \langle \od, \ox
\rangle.
\]
Hence
\[
|\langle d,x \rangle| \le \langle \od, \ox \rangle,
\]
and the last inequality follows from the definition of the norm of
an operator. \hfill $\square$

For bounding rigorously the optimal value, we claim boundedness
qualifications, which are more suitable for our purpose than
Slater's constraint qualifications. We assume that the conic
optimization problem satisfies the following condition which we
call {\it primal boundedness qualification} (PBQ):

\begin{itemize}
\item[(i)] Either the primal problem is infeasible,

\item[(ii)] or $\hf_p$ is finite, and there is a simple bound $\ox
\in \EK$ such that for every $\varepsilon > 0$ there exists a
primal feasible solution $x(\varepsilon)$ satisfying
$x(\varepsilon) \le \ox$ and $\langle c,x(\varepsilon)\rangle -
\hf_p \le \varepsilon$
\end{itemize}
Observe that PBQ implies that the primal problem is bounded from
below, but the existence of an optimal solution is not demanded,
only simple bounds $\ox$ for $\varepsilon$-optimal solutions are
required. The following theorem provides a finite lower bound
$\underline{f}_p$ of the primal optimal value.

\begin{theorem} \label{theo4.2}
Assume that PBQ holds. Let $\ty \in \EY^*$ and let $d:=
-A^*\ty+c$. Suppose further that $\ud^- \le \{d,0\}$, then:
\begin{itemize}
\item[(a)] The primal optimal value is bounded from below by
\begin{equation} \label{eq4.2}
\hf_p \ge \langle \ty, b \rangle + \langle \ud^-, \ox \rangle =:
\uf_p
\end{equation}

\item[(b)] If $\ud^- = 0$, then $\ty$ is dual feasible and $\hf_d
\ge \uf_p = \langle \ty,b\rangle$, and if moreover $\ty$ is
optimal then $\hf_d = \uf_p$.
\end{itemize}
\end{theorem}

{\sc Proof.} $(a)$ If the primal problem is infeasible, then
$\hf_p = + \infty$, and each finite value is a lower bound. Hence,
assume that PBQ (ii) is satisfied with $x:=x(\varepsilon)$ and
$\varepsilon>0$. Then
\[
\begin{array}{lcl}
\langle c,x \rangle & = & \langle d,x \rangle + \langle A^* \ty, x
\rangle \\
 & = & \langle \ty,b \rangle + \langle A^*\ty, x \rangle - \langle
 \ty,b \rangle + \langle d,x \rangle \\
 & = & \langle \ty, b \rangle + \langle \ty, Ax-b \rangle +
 \langle d,x \rangle.
\end{array}
\]
Since $x$ is primal feasible $Ax-b=0$, and
\[
\langle c,x \rangle = \langle \ty,b \rangle + \langle d,x \rangle.
\]
Lemma \ref{lem4.1} implies the inequality
\[
\langle c,x \rangle \ge \langle \ty,b \rangle + \langle \ud^-, \ox
\rangle.
\]
Because of PBQ (ii)
\[
\hf_p \ge \langle c,x \rangle - \varepsilon \ge \langle \ty,b
\rangle + \langle \ud^-, \ox \rangle- \varepsilon.
\]
For $\varepsilon \rightarrow 0$  the assertion (a) follows.

$(b)$ If $\ud^- = 0$ then $d \in \EK^*$, implying that $\ty$ is
dual feasible, and the assertion follows. \hfill $\square$

In particular, an approximate solution $\ty$ which is close to
optimality implies that $d$ is close to $\EK^*$. Hence, each lower
bound $\ud^-$ sufficiently close to $d^-$ is almost zero, and it
follows that $\uf_p \approx \langle \ty,b\rangle$ is reasonable;
that is the overestimation is not very much larger than necessary.
The lower bound uses the approximate optimal value $\langle \ty, b
\rangle$, and a correction is added which takes into account the
violation of dual feasibility $\ud^-$ evaluated at the upper bound
$\ox$.

We illustrate the bound for linear programming problems in
standard form
\begin{equation} \label{eq..}
\mbox{minimize } c^Tx \quad \mbox{s.t.} \quad Ax=b, \; x \ge 0.
\end{equation}
This is the special case of the conic optimization problem where
$\EX = \EX^* = \R^n$, $\EK = \EK^* = \R^n_+$ and $\EY = \EY^* =
\R^m$. It is well-known that for linear programming strong duality
$\hf_p = \hf_d =: \hf$ holds without any constraint
qualifications. Hence, Theorem \ref{theo4.2} yields immediately
the lower bound
\begin{equation} \label{eq*}
\hf \ge b^T\ty + (\ud^-)^T \ox =: \uf_p,
\end{equation}
where $\ud^-_j \le \min\{0, (-A^T \ty+c)_j\}$ for $j=1, \ldots,
n$. It is straightforward to control all effects of rounding
errors for computing $\uf_p$ by using directed rounding or
interval arithmetic.  The MATLAB toolbox  INTLAB \cite{Ru06b}
provides the directed rounding modes, and the following short
INTLAB program produces a rigorous lower bound:
\begin{verbatim}
     setround(-1);
     dlminus = min(0,A'*(-yt)+c);
     flow = b'*yt + dlminus'*xup;
     setround(0);
\end{verbatim}

If interval arithmetic is used, then the input data $A,b,c$ may be
intervals, and we obtain a lower bound for each instance within
the interval data. Verified error bounds for general linear
programming problems also with free variables can be found in
\cite{Ja04a}, and for formula (\ref{eq*}) see Corollary 6.1 in
\cite{Ja04a}.

To compute a rigorous upper bound of the optimal value we assume
that the conic optimization problem satisfies the following
condition, which we call the {\it dual boundedness qualification}
(DBQ):

\begin{itemize}
\item[(i)] Either the dual problem is infeasible,

\item[(ii)] or $\hf_d$ is finite, and there is a simple bound
$\oy$ such that for every $\varepsilon > 0$ there exists a dual
feasible solution $y(\varepsilon)$ satisfying $|y(\varepsilon)|
\le \oy$ and $\hf_d-\langle y(\varepsilon),b \rangle \le
\varepsilon$
\end{itemize}

\begin{theorem} \label{theo4.3}
Assume that DBQ holds. Let $\tx \in X$, and suppose further that
\begin{eqnarray}
|A\tx-b| & \le & \ore, \label{eq4.3}\\
\ux^- & \le & \{\tx, 0\}, \label{eq4.4} \; \mbox{and} \\
\od \; \; & \ge & -A^* y + c \label{eq4.5}\;  \mbox{for all dual
feasible} \; y \; \mbox{with} \; |y| \le \oy.
\end{eqnarray}
Then:
\begin{itemize}
\item[(a)] The dual optimal value is bounded from above by
\begin{equation} \label{eq4.6}
\hf_d \le \langle c,\tx \rangle - \langle \od, \ux^- \rangle +
\langle \oy , \ore \rangle =: \ohf_d.
\end{equation}

\item[(b)] If $\ux^- = 0$ and $\ore = 0$, then $\tx$ is primal
feasible and $\hf_p \le \ohf_d = \langle c, \tx \rangle$, and if
moreover $\tx$ is optimal, then $\hf_p = \ohf_d$.
\end{itemize}
\end{theorem}

{\sc Proof.} $(a)$ If the dual problem is infeasible then
$\hf_d=-\infty$, and each finite value is an upper bound. Hence,
assume that DBQ (ii) is satisfied with $y:=y(\varepsilon)$ and
$\varepsilon>0$. Then
\[
\begin{array}{lcl}
\langle y,b \rangle & = & -\langle y, A \tx-b \rangle + \langle
y,A\tx \rangle \\
 & = & \langle c,\tx \rangle + \langle y,A\tx \rangle - \langle
 c,\tx \rangle - \langle y, A\tx - b \rangle \\
 & = & \langle c,\tx \rangle - \langle c-A^*y, \tx \rangle -
 \langle y, A\tx-b \rangle.
\end{array}
\]
Since $y$ is dual feasible and $\od \ge d:= -A^*y+c \ge 0$, we can
apply Lemmata \ref{lem4.2} and \ref{lem4.3} which yield
\[
\langle y,b \rangle \le \langle c,\tx \rangle - \langle \od, \ux^-
\rangle + \langle \oy , \ore \rangle.
\]
Because of DBQ (ii)
\[
\hf_d \le \langle y,b \rangle + \varepsilon \le  \langle c,\tx
\rangle - \langle \od, \ux^- \rangle + \langle \oy , \ore  \rangle
+ \varepsilon.
\]
For $\varepsilon \rightarrow 0$ we obtain the upper bound
(\ref{eq4.6}).

The assertion $(b)$ follows immediately, since $\ux^- = 0$ and
 $\ore = 0$.   \hfill $\square$

Observe that for finite dimensional $y$ or in the case where $\EY$
is a vector lattice the absolute value $|.|$ is defined. If the
absolute value is not available, then we replace the inequalities
$|y(\varepsilon)| \le \oy, \quad |A\tx-b| \le \ore$ by
$\|y(\varepsilon)\| \le \oy$ and $\|A\tx-b\| \le \ore$,
respectively, and we obtain the error bound
\begin{equation} \label{eq4.9}
\hf_d \le \langle c,\tx \rangle - \langle \od,\ux^- \rangle +
 \oy \cdot \ore =: \of_d.
\end{equation}

Similarly as in the case of the lower bound, the upper bound uses
the approximate value $\langle c,\tx \rangle$ and takes into
account the violations of $\tx$ wrt. to the cone $\EK$ and the
linear equations.

The computation of the quantity $\langle \od,\ux^- \rangle$ can be
avoided. Since $\tx$ is an approximate optimal solution, $\tx \in
\EK$ or $\tx$ is close to $\EK$ (provided the conic solver has
computed reasonable approximations). Hence, for the supremum
$\tx^+ = \sup \{\tx,0\}$ the distance $\|\tx-\tx^+\|$ is small. If
we replace in Theorem \ref{theo4.3} $\tx$ by $\tx^+$ then the
quantity $|\langle c, \tx^+ \rangle - \langle c, \tx \rangle|$ is
small, but $\ux^-:=0 \le \{\tx^+,0\}$ yielding $\langle \od,\ux^-
\rangle =0$ and the upper bound
\begin{equation} \label{eq31}
\hf_d \le \langle c, \tx^+ \rangle + \langle \oy, \ore \rangle,
\end{equation}
where $|A\tx^+ -b| \le \ore$. In general it is not possible to
compute the supremum $\tx^+$ exactly, but each close upper bound
$\stackrel{\approx}{x} \ge \tx^+$ will suffice.

In the special case of linear programming we can take the exact
supremum $\tx^+ \in \R^{n}_+$ defined by (\ref{eq2.8}). Then
$\ux^- = 0$, and we obtain the upper bound
\begin{equation} \label{eq4.9l}
\hf \le c^T \tx^+  + \oy^T \ore = \of_d.
\end{equation}
The following short INTLAB program produces this upper bound:

\begin{verbatim}
     xtplus = max(0,xt)
     setround(-1);
     rn = abs(A*xtplus -b);
     setround(+1);
     rp = abs(A*xtplus -b);
     r = max(rn,rp);
     fu = c'*xtplus + yup'*r;
     setround(0);
\end{verbatim}

Until now we have assumed the existence of $\varepsilon$-optimal
solutions within some reasonable bounds. Now we mention briefly
that in the case where appropriate primal or dual boundedness
qualifications are not known it is frequently possible to compute
verified primal and dual feasible solutions which are close to
optimality. These solutions can be used to compute verified
reasonable error bounds for the optimal value. The basic algorithm
consists of the following steps:

\begin{itemize}
\item[(i)] Perturb the original problem slightly such that the
optimal solution of the perturbed problem is an interior feasible
solution of the original problem.

\item[(ii)] Solve the perturbed problem approximately.

\item[(iii)] Use this approximation to compute an enclosure (i.e
an appropriate interval) containing a feasible solution.

\item[(iv)] Evaluate the objective function for the enclosure.
\end{itemize}

Especially step (iii) is nontrivial since the existence of
feasible solutions must be rigorously proved. Interval arithmetic
provides several methods for computing enclosures of solutions for
linear and nonlinear systems in the finite dimensional case, but
also for infinite dimensional problems (certain types of ordinary
and partial differential equations) enclosure methods are known.
However the bounds obtained in this way have two disadvantages.
They are much more time-consuming than the previous ones, and they
provide an upper bound of the primal optimal value only if the
primal problem is well-posed, and a lower bound of the dual
optimal value only if the dual problem is well-posed.

A detailed description of this algorithm can be found in the case
of linear programming in \cite{Ja04a}, for smooth convex
programming problems in \cite{Ja04}, and for semidefinite
programming problems and linear matrix inequalities in
\cite{JaChaKeil05}.

\section{Second Order Cone Programming}

In SOCP the partial ordering is defined by the ice-cream cone
(\ref{eq2.9}) yielding a finite dimensional vector lattice
equipped with the following operations.

\begin{theorem} \label{theo5.1}
Let $\EK \subseteq \R^n$ be defined by (\ref{eq2.9}). Then for $x
\in \R^n$
\begin{equation} \label{eq5.1*}
x^+= \left\{ \begin{array}{ll} x & \mbox{if} \quad x_n \ge \|x_:\|\\[0.6cm]
0  & \mbox{if} \quad x_n \le -\|x_:\|\\[0.6cm]
\dfrac{\|x_{:}\|+x_n}{2\|x_:\|}{\displaystyle {x_: \choose
\|x_:\|}} & \mbox{if} \quad -\|x_:\|< x_n < \|x_:\|
\end{array} \right.
\end{equation}
and
\begin{equation} \label{eq5.2*}
x^-= \left\{ \begin{array}{ll} x & \mbox{if} \quad x_n \le -\|x_:\|\\[0.6cm]
0  & \mbox{if} \quad x_n \ge \|x_:\|\\[0.6cm]
-\dfrac{\|x_{:}\|-x_n}{2\|x_:\|}{\displaystyle {-x_: \choose
\|x_:\|}} & \mbox{if} \quad - \|x_:\| < x_n < \|x_:\|
\end{array} \right.
\end{equation}
\end{theorem}

{\sc Proof.} First we prove (\ref{eq5.1*}). If $x_n \ge \|x_:\|$
then $x  \in \EK$ which implies $x^+:= \sup\{x,0\}=x$. If $x_n \le
-\|x_:\|$ then $x \in -\EK$ and $x^+ = 0$. Finally, assume that
$-\|x_:\| < x_n < \|x_:\|$. Then $x_: \neq 0$, and a simple
geometric argument shows that $x^+$ is the orthogonal projection
of $x$ onto the boundary of $\EK$, that is
\[
x^+=\alpha {x_: \choose \|x_:\|} \quad \mbox{and} \quad
0=(x^+-x)^T x^+.
\]
The latter condition describes the orthogonality. Since
\[
\begin{array}{lcl} 0 & = & \left( \begin{array}{ll} \alpha x_: & -x_: \\ \alpha
\|x_:\| & -x_n
\end{array}\right)^T \left( \begin{array}{l}
\alpha  x_: \\ \alpha  \|x_:\|
\end{array} \right)\\[0.5cm]
 & = & \alpha^2 \|x_:\|^2-\alpha\|x_:\|^2 + \alpha^2 \|x_:\|^2 -
 \alpha x_n \|x_:\|\\[0.5cm]
 & = & 2 \alpha^2 \|x_:\|^2 - \alpha(\|x_:\|^2+x_n \|x_:\|),
 \end{array}
\]
we get $\alpha = (\|x_:\|+x_n)/2\|x_:\|$ which proves
(\ref{eq5.1*}).

Finally, (\ref{eq5.2*}) follows from (\ref{eq5.1*}) since $x^- =
\inf\{x,0\}=-\sup\{-x,0\}$. \hfill $\square$

Due to rounding-off errors $x^-$ may not be computed exactly. But
for computing rigorous results we know from the previous section
that it is sufficient to compute a lower bound $\ux^- \le x^-$ by
using directed rounding or interval arithmetic. The distinction of
cases, however, must be implemented carefully when $x_n$ is almost
equal to $-\|x_:\|$ or $\|x_:\|$. A similar remark applies to the
computation of an upper bound $\ox^+ \ge x^+$.

Let $\EK$ be the Cartesian product of ice-cream cones $\EK_j
\subseteq \R^{n_j}$ for $j=1, \ldots, n$. This is a convex,
self-dual cone (see Alizadeh, Goldfarb \cite{AlGo03}). The
\emph{standard SOCP problem} has the form
\begin{equation} \label{eq5.3*}
\mbox{minimize } \sum\limits_{j=1}^n c_j^T x_j \quad \mbox{s.t. }
\sum\limits_{j=1}^n A_jx_j = b, \; x_j \in \EK_j \quad \mbox{for }
j=1, \ldots, n,
\end{equation}
where $A_j \in \R^{m \times n_j}$, $c_j, x_j \in \R^{n_j}$ and $b
\in \R^m$. If we merge these quantities
\begin{equation} \label{eq5.3a*}
\begin{array}{lcl}
A & := & (A_1; \ldots; A_n),\\[0.3cm]
c & := & (c_1; \ldots; c_n),\\[0.3cm]
x & := & (x_1; \ldots; x_n),\\
\end{array}
\end{equation}
then the standard SOCP problem has the form (\ref{eq3.1}), and it
follows that the dual problem (\ref{eq3.5}) can be written as
\begin{equation} \label{eq5.4*}
\mbox{minimize } b^Ty \quad \mbox{s.t. } -A_j^T y+c_j \in \EK_j
\quad \mbox{for } j=1, \ldots, n.
\end{equation}
Here we have chosen the finite-dimensional spaces $\EX:=\R^{\on}$
where $\on = \Sigma_jn_j$ and $\EY:= \R^m$ equipped with the
Euclidean inner products. By $x_{i,j}$ we denote the $i$-th
component of the vector $x_j$, and $x_{:,j}:= (x_{1,j} \ldots,
x_{n_j-1,j})^T$. In this section
\[
\ox = (\ox_1; \ldots; \ox_n)
\]
denotes a vector in $\EK$ with $\ox_{:,j}=0$ and $\ox_{n_j,j}
> 0$ for every $j$. Then
\begin{equation} \label{eq5.5*}
x \le \ox \quad \Leftrightarrow \quad \|x_{:,j}\|+x_{n_j,j} \le
\ox_{n_j,j} \quad \mbox{for } j=1, \ldots, n.
\end{equation}
The computation of a rigorous lower bound for the optimal value of
(\ref{eq5.3*}) is a straightforward application of Theorem
\ref{theo4.2}.
\begin{coro} \label{coro5.1}
Assume that PBQ holds for some $\ox \in \EK$. Let $\ty \in \R^m$,
and let
\begin{equation} \label{eq5.6*}
d_j:=-A_j^T \ty +c_j \quad \mbox{for} \; j=1, \ldots, n.
\end{equation}
Suppose further that for $j=1, \ldots, n$ there are lower bounds
$\ud_{j}^- \le d_{j}^-$. Then:
\begin{itemize}
\item[(a)] The primal optimal value is bounded from below by
\begin{equation} \label{eq5.7*}
\hf_p \ge b^T\ty + \sum\limits_{j=1}^n \ud^-_{n_j,j} \ox_{n_j,j}:=
\uf_p.
\end{equation}

\item[(b)] If $d_{n_j,j} \ge \|d_{:,j}\|$ for $j = 1, \ldots, n$,
then $\ty$ is dual feasible and $\hf_d \ge \uf_p = b^T\ty$, and if
moreover $\ty$ is optimal then $\hf_d = \uf_p$.
\end{itemize}
\end{coro}

{\sc Proof.}

 It follows from (\ref{eq4.2}) that
\[
\begin{array}{lcl}
\hf_p & \ge & \langle \ty,b \rangle + \langle \ud^-, \ox \rangle\\
 & = & b^T\ty + \sum\limits_{j=1}^n \langle \ud_j^-, \ox_j \rangle
 \\
  & = & b^T\ty + \sum\limits_{j=1}^n \ud^-_{n_j,j} \ox_{n_j, j},
\end{array}
\]
where the last equation is fulfilled because $\ox_{:,j}=0$. This
finishes the proof of $(a)$. If $d_{n_j,j} \ge \|d_{:,j}\|$ for $j
= 1, \ldots, n$ then $-A_j^T \ty + c_j \in \EK_j$, $d_j^- = 0$,
and $\ty$ is dual feasible. Hence, $(b)$ is proved. \hfill
$\square$

An upper bound for the optimal value of (\ref{eq5.3*}) is an
immediate application of Theorem \ref{theo4.3}.

\begin{coro} \label{coro5.2}
Assume that DBQ is fulfilled. Let $\tx \in \EK$, and suppose
further that
\[
|\sum\limits_{j=1}^n A_j \tx_j-b| \le \ore,
\]
then:
\begin{itemize}
\item[(a)] The dual optimal value is bounded from above by
\[
\hf_d \le \sum\limits_{j=1}^n c_j \tx_j + \oy^T\ore=: \of_d.
\]

\item[(b)] If $\ore =0$ then $\tx$ is primal feasible and $\hf_p
\le \of_d = \sum_{j=1}^n c_j \tx_j$, and if moreover $\tx$ is
optimal then $\hf_p = \of_d$.
\end{itemize}
\end{coro}

{\sc Proof.} Since $\tx \in \EK$, it follows that $\ux^- = 0 \le
\{\tx,0\}$. Therefore the assertion is an immediate consequence of
Theorem \ref{theo4.3}. \hfill $\square$

SOCP solvers may compute an approximation $\tx$ which is not
contained in $\EK$, i.e. for some $j$ the part $\tx_j$ is not
contained in the convex cone $\EK_j$. Then $\tx$ is replaced by an
upper bound of $\tx^+$. As aforementioned, in floating point
arithmetic formula (\ref{eq5.1*}) must be carefully evaluated
using directed rounding such that the computed result is
guaranteed to be in $\EK$.

\section{Semidefinite Programming}

We examine the \emph{standard primal semidefinite programming
problem}
\begin{equation} \label{eq6.1}
\mbox{minimize } \langle C,X \rangle \quad \mbox{s.t. } \langle
A_i,X\rangle = b_i, \quad i=1, \ldots, m, \quad X \in \EK,
\end{equation}
where $C,X$ and $A_i$ are real symmetric $s \times s$ matrices, $b
\in \R^m$, $\EK = S_+^s$, and $\langle \;,\rangle$  denotes the
inner product (\ref{eqz2}) in the linear space of symmetric
matrices. Using the svec operator (\ref{eqz1}) such that
\begin{equation} \label{eq6.2}
c=\svec (C), \quad x = \svec(X), \quad a_i = \svec(A_i),
\end{equation}
we can write problem (\ref{eq6.1}) equivalently in the form
\begin{equation} \label{eq6.3}
\mbox{minimize } c^Tx \quad \mbox{s.t. } Ax=b, \quad x \in \EK,
\end{equation}
where $A$ is the matrix with rows $a_i^T$. Problem (\ref{eq6.3})
has the form (\ref{eq3.1}), and therefore the dual problem
(\ref{eq3.5}) is
\begin{equation} \label{eq6.4}
\mbox{maximize } b^Ty \quad \mbox{s.t. } -A^Ty+c \in \EK, \quad y
\in \R^m.
\end{equation}
The equivalent matrix notation is
\begin{equation} \label{eq6.5}
\mbox{maximize } b^Ty \quad \mbox{s.t. } -\sum\limits_{i=1}^m y_i
A_i + C \in \EK.
\end{equation}

\begin{coro} \label{coro6.1}
Assume that the SDP satisfies:
\begin{itemize}
\item[(i)] Either the primal problem is infeasible, or

\item[(ii)] there exists a nonnegative number $\ox$ such that for
every $\varepsilon > 0$ there exists a primal feasible solution
$X(\varepsilon) \le \ox \cdot I$ and $\langle
C,X(\varepsilon)\rangle -\hf_p \le \varepsilon$,
\end{itemize}
where I denotes the identity matrix. Let $\ty \in \R^m$, and let
\begin{equation} \label{eq6.5*}
D=C-\sum\limits_{i=1}^m \ty_i A_i.
\end{equation}
Suppose further that $\ud^- \le \{\lambda_{\min}(D),0\}$, and that
$D$ has at most $l$ negative eigenvalues. Then:
\begin{itemize}
\item[(a)] The primal optimal value is bounded from below by
\begin{equation} \label{eq6.6}
\hf_p \ge b^T \ty+ l \cdot \ud^- \cdot \ox =: \uf_p.
\end{equation}
\item[(b)] If $\ud^- = 0$ then $\ty$ is dual feasible and $\hf_d
\ge \uf_p = b^T \ty$, and if moreover $\ty$ is optimal then $\hf_d
= \uf_p$.
\end{itemize}
\end{coro}

{\sc Proof.} Let $D$ have the eigenvalue decomposition $D=Q^T
\Lambda Q$, then (\ref{eq2.12}) implies $D^-=Q^T \Lambda^- Q$.
Hence $\uD^-:= \ud^- \cdot Q^T \mbox{sign}(-\Lambda^-) Q \le
\{D,0\}$, where sign  is $+1, 0$ or $-1$ if the corresponding
coefficient of the matrix is positive, zero or negative,
respectively. Moreover, PBQ implies $\oX:=\ox \cdot I \ge
X(\varepsilon)$. Now from Theorem \ref{theo4.2} (a) it follows
that
\[
\hf_p \ge \langle \ty, b\rangle + \langle \uD^-, \oX \rangle = b^T
\ty + l \cdot \ud^- \cdot \ox
\]
which proves (a). The part (b) is an immediate consequence of
Theorem \ref{theo4.2} (b). \hfill $\square$

In order to control all rounding errors and to compute a verified
lower bound $\uf_p$ it is necessary to compute a rigorous lower
bound of the smallest eigenvalue for a symmetric matrix.
Interesting references for computing rigorous bounds of some or
all eigenvalues are Floudas \cite{Fl2000}, Mayer \cite{May94a},
Neumaier \cite{Neumaier9}, and Rump \cite{Ru93,Ru94}. In VSDP  we
have computed the quantities $l$ and $\ud^-$ by using Weyl's
Perturbation Theorem for symmetric matrices: For an approximate
eigenvalue decomposition $\tD=\tQ^T \tLambda \tQ$ of $D$ with
eigenvalues $\tlambda_i$ on the diagonal of $\tLambda$, we use
directed rounding or interval arithmetic for computing an error
matrix $E \ge |D-\tD|$. Then the Theorem of Weyl implies that
\[
|\lambda_i(D)-\tlambda_i| \le \|E\|_2
\]
for each eigenvalue $\lambda_i(D)$. Therefore, we obtain the
bounds
\begin{equation} \label{eq6.6a}
\tlambda_i - \|E\|_{\infty} \le \lambda_i(D) \le
\tlambda_i+\|E\|_{\infty}
\end{equation}
for all eigenvalues, yielding immediately the quantities $\ud^-$
and $l$. The short INTLAB program

\begin{verbatim}
     [Qt,Lt] = eig(full(mid(D)));
     E = D - Qt * intval(Lt) * Qt';
     r = abss(norm(E,inf));
     lambda = midrad(diag(Lt),r);
\end{verbatim}

implements these bounds for the eigenvalues of a symmetric
interval matrix {\tt D}, where {\tt eig}, {\tt abss} and {\tt
midrad} denote the MATLAB and INTLAB routines for computing
approximate eigenvalues, the absolute value, and the midpoint
radius description of interval quantities, respectively.

\begin{coro} \label{coro6.2}
Assume that DBQ is fulfilled. Let $\tX \in \EK = S_+^s$ and
suppose further that
\begin{equation} \label{eq6.7}
|\langle A_i, \tX \rangle-b_i| \le \ore_i \quad \mbox{for } i=1,
\ldots, m.
\end{equation}
Then:
\begin{itemize}
\item[(a)] The dual optimal value is bounded from above by
\begin{equation} \label{eq6.8}
\hf_d \le \langle C,\tX \rangle + \oy^T \ore =: \of_d.
\end{equation}
\item[(b)] If $\ore = 0$ then $\tX$ is primal feasible and $\hf_p
\le \of_d = \langle C,\tX \rangle$, and if moreover $\tX$ is
optimal then $\hf_p = \of_d$.
\end{itemize}
\end{coro}

{\sc Proof.} Because $\tX^-=0$, this corollary is an immediate
consequence of Theorem \ref{theo4.3}. \hfill $\square$

In general, SDP-solvers, compute an approximation $\tX$ which is
not a positive semidefinite matrix, but there is a cluster of
negative eigenvalues of $\tX$ nearby zero. In order to enforce a
positive  semidefinite approximation, we compute a rigorous bound
$\ux \le \{\lambda_{\min}(\tX),0\}$. Then it follows that $\tX-\ux
I \in \EK$, i.e. is positive semidefinite, and we can use this
shifted matrix in the previous corollary. Then with directed
rounding it is straightforward to compute $\of_d$.

\section{Block Structured Variables}

Frequently conic optimization problems have block structured
variables, that is the variables are in the Cartesian product of
different cones. More precisely, there are $n$ real normed vector
spaces $\EX_1, \ldots, \EX_n$, convex cones $\EK_1 \subseteq
\EX_1, \ldots, \EK_n \subseteq \EX_n$, a real normed vector space
$\EY$, and $n$ continuous linear operators $A_j : \EX_j
\rightarrow \EY$. Let $\EX$ and $\EK$ denote the Cartesian
products of the spaces $\EX_j$ and the cones $\EK_j$,
respectively. The vectors $x$ and $c$ and the linear operator $A$
are partitioned appropriately:
\[
\begin{array}{lclll}
x & = & (x_1; \ldots; x_n) & \mbox{where} & x_j \in \EX_j,\\[0.3cm]
c & = & (c_1; \ldots; c_n) & \mbox{where} & c_j \in \EX^*_j,\\[0.3cm]
A & = & (A_1; \ldots; A_n). & &
\end{array}
\]
Defining
\begin{equation} \label{eq7.0}
Ax := \sum\limits_{j=1}^n A_j x_j \quad \mbox{and} \quad \langle
c,x \rangle := \sum\limits_{j=1}^n \langle c_j, x_j \rangle,
\end{equation}
it follows that $A : \EX \rightarrow \EY$ is a continuous linear
operator, and $c \in \EX^*$. The {\emph {primal conic optimization
problem with block structured variables}} has the form
\begin{equation} \label{eq7.1}
\mbox{minimize } \sum\limits_{j=1}^n \langle c_j, x_j \rangle
\quad \mbox{s.t. } \sum\limits_{j=1}^n A_j x_j = b, \quad x_j \in
\EK_j \quad \mbox{for} \; j=1, \ldots, n,
\end{equation}
and the dual problem is
\begin{equation} \label{eq7.2}
\mbox{minimize }  \langle y,b \rangle \quad \mbox{s.t. } (-A_1^*
y; \ldots; -A_n^*y)+(c_1; \ldots ; c_n) \in \EK_1^* \times \ldots
\times \EK_n^*, \quad y \in \EY^*.
\end{equation}
The most important examples are
\emph{semidefinite-quadratic-linear programs}. These are block
structured problems where $\EY = \R^m$,
\begin{equation} \label{eq7.3}
x = (x_1^s; \ldots; x_{n_s}^s; x_1^q; \ldots; x_{n_q}^q; x^l) \in
\EX,
\end{equation}
and
\[
\begin{array}{ll}
x_1^s, \ldots, x_{n_s}^s & \mbox{are symmetric matrices of various
sizes,}\\[0.3cm]
x_1^q, \ldots, x_{n_q}^q & \mbox{are vectors of various
sizes, and}\\[0.3cm]
x^l & \mbox{is a vector.}
\end{array}
\]
Using (\ref{eq7.0}), (\ref{eq7.1}), (\ref{eq5.3*}) and
(\ref{eq6.1}) we obtain the primal problem
\begin{equation} \label{eq7.4}
\begin{array}{rl}
\mbox{minimize} & \sum\limits_{j=1}^{n_s} \langle c_j^s, x_j^s
\rangle  +  \sum\limits_{k=1}^{n_q}(c_k^q)^Tx_k^q +
(c^l)^Tx_l\\[0.5cm]
\mbox{s.t.} & \sum\limits_{j=1}^{n_s} \langle A_{ij}^s, x_j^s
\rangle  + \sum\limits_{k=1}^{n_q} (A_{ik}^q)^T x_k^q +
(A^l_i)^T x_l = b_i, \quad i=l, \ldots, m\\[0.5cm]
 & \multicolumn{1}{c}{x_j^s \in \EK_j^s, \quad x_k^q \in \EK_k^q, \quad x^l \in \EK^l \quad \forall \; j,k,}
\end{array}
\end{equation}
where the matrices and vectors have appropriate dimensions, and
$\EK_j^s$, $\EK^q_k$ and $\EK^l$ are the convex cones of positive
semidefinite matrices, the ice-cream cones, and the positive
orthant, respectively.

The dual problem has the form
\begin{equation} \label{eq7.5}
\begin{array}{rll}
\mbox{maximize} & \sum\limits_{i=1}^m b_iy_i & \\[0.5cm]
\mbox{s.t.} & -\sum\limits_{i=1}^m A_{ij}^s y_i + c_j^s \in \EK_j^s & \mbox{for } j=1, \ldots, n^s \\[0.5cm]
 & -\sum\limits_{i=1}^m A_{ik}^q y_i + c_k^q \in \EK_k^q & \mbox{for } k=1, \ldots, n^q \\[0.5cm]
  & -\sum\limits_{i=1}^m A_{i}^l y_i + c^l \in \EK^l. &  \\[0.5cm]
\end{array}
\end{equation}
Observe that the set of primal feasible solutions is the Cartesian
product of semidefinite, quadratic and nonnegative orthant cones
intersected with an affine subspace. It is possible that $n^s$,
$n^q$ or the length of $x^l$ is zero, which means that one or more
of the three parts of the problem is absent.

The following two corollaries provide finite lower and upper
bounds of the optimal value for block structured problems.

\begin{coro} \label{coro7.1}
Assume that PBQ holds for some $\ox = (\ox_1; \ldots; \ox_n) \in
\EK$. Let $\ty \in \EY$, and assume that for $j=1, \ldots, n$
\begin{equation} \label{eq7.6}
d_j:= -A_j^* \ty + c_j
\end{equation}
and $\ud_j^- \le \{d_j, 0\}$. Then:
\begin{itemize}
\item[(a)] The primal optimal value is bounded from below by
\begin{equation} \label{eq7.7}
\hf_p \ge \langle \ty, b \rangle + \sum\limits_{j=1}^n \langle
\ud_j^-, \ox_j \rangle =: \uf_p.
\end{equation}
\item[(b)] If $\ud_j^- = 0$ for $j = 1, \ldots, n$, then $\ty$ is
dual feasible and $\hf_d \ge \uf_p = \langle \ty, b \rangle$, and
if moreover $\ty$ is optimal then $\hf_d = \uf_p$.
\end{itemize}
\end{coro}
{\sc Proof.} This corollary follows immediately from Theorem
\ref{theo4.2} by observing the linearity of the block structured
variables. \hfill $\square$

\begin{coro} \label{coro7.2}
Assume that DBQ holds for some $\oy \in \EY^*$. Let $\tx = (\tx_1;
\ldots; \tx_n) \in \EK$ and suppose further that
\begin{equation} \label{eq7.8}
|\sum\limits_{j=1}^n A_j \tx_j - b| \le \ore,
\end{equation}
then
\begin{itemize}
\item[(a)] The dual optimal value is bounded from above by
\begin{equation} \label{eq7.8*}
\hf_d \le \sum\limits_{j=1}^n \langle c_j, \tx_j \rangle + \langle
\oy, \ore \rangle =: \of_d.
\end{equation}
\item[(b)]  If $\ore = 0$ then $\tx$ is primal feasible and $\hf_p
\le \of_d$, and if moreover $\tx$ is optimal then $\hf_p = \of_d$.
\end{itemize}
\end{coro}
{\sc Proof.} The corollary is an immediate consequence of Theorem
\ref{theo4.3}. \hfill $\square$

Lower and upper bounds for the optimal value of
semidefinite-quadratic-linear programs can be immediately obtained
by inserting the preciding formulas for LP, SOCP and SDP.

\section{Certificates of Infeasibility}

A theorem of alternatives states that for two systems of equations
or inequalities, one or the other system has a solution, but not
both. A solution of one of the systems is called a
\emph{certificate of infeasibility} for the other which has no
solution, since in principle this allows an easy check to prove
infeasibility. Certificates of infeasibility are frequently
computed by optimization algorithms if no feasible solutions of
the primal or dual constraints exist. Especially in the presence
of equality constraints, certificates cannot be represented
exactly in floating point arithmetic, and therefore approximate
certificates can satisfy the constraints only within certain
tolerances. This effect is amplified by rounding-off errors during
the calculations for computing the approximate certificate.
However, it turns out that in order to prove infeasibility by
using floating-point arithmetic it is sufficient if an interval of
small diameter can be given which guarantees to contain a
certificate. We call such an interval a \emph{rigorous} or
\emph{verified certificate of infeasibility}, and describe briefly
how such certificates can be obtained for conic problems. We begin
with two well known propositions and include the short proofs.

\begin{prop} \label{prop5.1}
Suppose that $\ty \in Y^*$ satisfies
\begin{equation} \label{eq5.1}
A^* \ty \in \EK^*, \quad \langle \ty,b \rangle <0,
\end{equation}
then the system of primal constraints
\begin{equation} \label{eq5.2}
Ax=b, \quad x \in \EK
\end{equation}
has no solution.
\end{prop}

{\sc Proof.} If the system (\ref{eq5.2}) has a solution $x$, then
$0 \le \langle A^*\ty, x \rangle = \langle \ty,Ax \rangle =
\langle \ty, b \rangle$ contradicting our assumption $\langle \ty,
b \rangle <0$. \hfill $\square$

The linear functional $\ty$ is called a \emph{certificate} of
\emph{primal infeasibility}, and represents a dual unbounded ray.

\begin{prop} \label{prop5.2}
Suppose that $\tx \in X$ satisfies
\begin{equation} \label{eq5.3}
A\tx = 0, \quad \tx \in \EK, \quad \langle c,\tx \rangle <0
\end{equation}
then the system of dual constraints
\begin{equation} \label{eq5.4}
-A^*y+c \in \EK^*, \; y \in Y^*
\end{equation}
has no solution.
\end{prop}

{\sc Proof.} If the system (\ref{eq5.4}) has a solution $y \in
Y^*$, then $0 \le \langle -A^*y+c, \tx \rangle = - \langle y,A\tx
\rangle + \langle c,\tx \rangle = \langle c,\tx \rangle <0$
contradicting our assumption. Hence, system (\ref{eq5.4}) has no
solution. \hfill $\square$

The vector $\tx$ is called a \emph{certificate} of \emph{dual
infeasibility} and represents a primal unbounded ray.

Many conic solvers expose infeasibility by computing approximate
unbounded rays. Given an approximate primal unbounded ray $\tx \in
\EX$, dual infeasibility is proved if the equation and sign
conditions (\ref{eq5.3}) can be checked rigorously on a computer.
The underdetermined equation $A \tx = 0$ is in general not exactly
satisfied for floating-point certificates $\tx$. To obtain a
rigorous certificate we proceed as follows: Let $\beta$ be an
approximation of $\langle c, \tx \rangle$, and assume that $\beta
< 0$. Otherwise, for nonnegative $\beta$ the sign condition in
(\ref{eq5.3}) would be even not satisfied for the approximate
primal unbounded ray $\tx$, and this indicates that $\tx$ is not
suitable. Then we compute an interval of small diameter {\bf x},
also called \emph {enclosure}, for a solution of the
underdetermined linear system
\begin{equation} \label{eq38}
Ax = 0 \quad \mbox{and} \quad \langle c,x \rangle = \beta < 0,
\end{equation}
which is close to $\tx$. If $\mbox{{\bf x}} \subseteq \EK$, then
there exists an $\stackrel{\approx}{x} \; \in {\bf x}$ which
fulfills the condition (\ref{eq5.3}) yielding a rigorous
certificate of dual infeasibility. This check depends on the
special problem and requires further information about the
operator $A$ and the cone $\EK$. In the three cases LP, SOCP, and
SDP  for  the finite-dimensional linear system (\ref{eq38})
methods of interval arithmetic can be used for computing an
appropriate enclosure {\bf x}. For a detailed description of such
an algorithm see \cite{Ja04a}. The condition $\mbox{{\bf
x}}=[\ux,\ox] \subseteq \EK$ can be verfied for LP by checking the
equivalent condition
\begin{equation} \label{eq35a}
\ux \ge 0,
\end{equation}
for SOCP we check the equivalent condition
\begin{equation} \label{eq36a}
\underline{x}_n \ge \|\underline{x}_{:}\|,
\end{equation}
and for SDP we check
\begin{equation} \label{eq37a}
\lambda_{\min}(\underline{X}) \ge 0,
\end{equation}
where $\underline{x}= \svec(\underline{X})$.

In the case of an approximate dual improving ray $\ty$, primal
infeasibility can be rigorously proved on a computer if the
condition (\ref{eq5.1}) can be verified; that is, if the sign
conditions
\begin{equation} \label{eq34}
\langle \ty, b \rangle < 0 \quad \mbox{and} \quad \langle A^* \ty,
x \rangle \ge 0 \quad \mbox{for all } x \in \EK
\end{equation}
can be checked reliably. As before, this check depends on the
special problem and requires further information about the
operator $A^*$ and the cone $\EK^*$.

It follows immediately that for LP (\ref{eq34}) is equivalent to
\begin{equation} \label{eq35}
b^T \ty < 0 \quad \mbox{and} \quad A^T \ty \ge 0,
\end{equation}
for SOCP we obtain the equivalent condition
\begin{equation} \label{eq36}
b^T \ty < 0 \quad \mbox{and} \quad (A_j^T \ty)_{n_j} \ge \|(A_j^T
\ty)_:\| \quad \mbox{for } j = 1, \ldots, n,
\end{equation}
and for SDP we get
\begin{equation} \label{eq37}
b^T \ty <0 \quad \mbox{and} \quad
\lambda_{\min}(\sum\limits_{i=1}^m \ty_i A_i) \ge 0.
\end{equation}
All three conditions can be checked rigorously by using directed
rounding or interval arithmetic. The vector $\ty$ provides a
rigorous certificate which can be viewed as a degenerate interval
with zero diameter.

\section{Combinatorial Optimization}

Linear and semidefinite programs play a very useful role in global
and combinatorial optimization (Wolkowicz \cite{Wolko}). Several
methods (for example lift-and-project methods) are known for
constructing linear or semidefinite relaxations, which are used in
branch-bound-and-cut algorithms to eliminate regions which do not
contain global minimizers. Neumaier and Shcherbina
\cite{Neumaier5} have pointed out that backward error analysis has
no relevance for combinatorial programs, since slightly perturbed
coefficients no longer produce problems of the same class. There,
one can also find an innocent-looking linear integer problem for
which the commercial high quality solver CPLEX \cite{CPLEX} and
several other state-of-the-art solvers failed. The reason is that
the relaxations are not solved with sufficient accuracy and global
minimizers are truncated. Hence, in order to obtain safe results,
it is important to have reliable, good and cheaply computable
lower bounds of the optimal value for relaxations.

Various problems like Max-Cut, Partitioning, Coloring and many
others can be formulated as linear integer problems where the
vector of decision variables $x \in \{-1,1\}^n$. Tight
semidefinite relaxations are obtained by lifting the vector $x$
into the space of semidefinite matrices by the operation
\begin{equation} \label{eq91}
X = xx^T.
\end{equation}
It follows immediately that
\begin{equation} \label{eq92}
X \succeq 0, \; \diag(X) = e, \; \mbox{and} \; \rank(X) = 1,
\end{equation}
where $e$ is the vector of ones. Dropping the condition $\rank(X)
= 1$ we obtain a semidefinite relaxation. Laurent and Poljak
\cite{LauPol95} have shown that for this type of relaxations $-1
\le X_{ij} \le 1$, and if $X_{ij} \in \{-1,1\}$ then $X = xx^T$
where $x \in \{-1,1\}^n$. This property establishes the tightness
of these relaxations. Moreover, it follows that the primal
boundedness qualification is fulfilled in the way that an optimal
solution exists with $\lambda_{\max}(X) \le n$, and thus  rigorous
lower bounds for the optimal value can be computed.

Sometimes, these tight relaxations are in addition ill-posed. As
an example we consider Graph Partitioning Problems. These are
known to be NP-hard, and finding an optimal solution is difficult.
Graph Partitioning has many applications among those is VLSI
design. Here, we investigate semidefinite relaxations for the
special case of Equicut Problems, which have turned out to deliver
tight lower bounds (see also Gruber and Rendl \cite{GruRe02}). The
general case of Graph Partitioning Problems can be treated
similarly.

Given an edge-weighted graph $G$ with an even number of vertices,
the problem is to find a partitioning of the vertices into two
sets of equal cardinality which minimizes the weight of the edges
joining the two sets. The algebraic formulation is obtained by
representing the partitioning as an integer vector $x \in \{-1,
1\}^n$ satisfying the parity condition $\sum_i x_i=0$. Then the
Equicut Problem is equivalent to
\[
\min\sum\limits_{i<j} a_{ij} \dfrac{1-x_ix_j}{2} \quad
\mbox{subject to} \quad x \in \{-1,1\}^n, \; \sum\limits_{i=1}^n
x_i = 0,
\]
where $A=(a_{ij})$ is the symmetric matrix of edge weights. This
follows immediately, since $1-x_ix_j=0$ iff the vertices $i$ and
$j$ are in the same set. The objective can be written as
\[
\dfrac{1}{2} \sum\limits_{i<j} a_{ij}(1-x_ix_j)= \dfrac{1}{4}
x^T(\Diag(Ae)-A)x = \dfrac{1}{4} x^TLx,
\]
where $L:= \Diag(Ae)-A$ is the {\it Laplace matrix} of $G$. Using
$x^TLx=\trace(L(xx^T))$ and $X=xx^T$, it can be shown that this
problem is equivalent to
\[
\hf_p = \min \dfrac{1}{4} \langle L,X \rangle \mbox{ subject to }
 \diag(X) = e, \; e^TXe = 0, \; X \succeq 0, \; \rank(X)=1.
\]
Since $X\succeq 0$ and $e^TXe=0$ implies $X$ to be singular, the
problem is ill-posed, and for arbitrarily small perturbations of
the right hand side the problem becomes infeasible. By definition,
the Equicut Problem has a finite optimal value $\hf_p$, and a
rigorous upper bound of $\hf_p$ is simply obtained by evaluating
the objective function for a given partitioning integer vector
$x$. Hence, it is left over to compute a rigorous lower bound. At
first, the nonlinear rank one constraint is left out yielding  an
ill-posed semidefinite relaxation, where the Slater condition does
not hold. The related constraints $\diag(X)=e$ and $e^TXe=0$ can
be written as
\[
\langle A_i, X \rangle = b_i, \; b_i=1, \;A_i=E_i \; \mbox{for} \;
i=1, \ldots, n, \; \mbox{and} \; A_{n+1} = ee^T, \; b_{n+1} = 0.
\]
where $E_i$ is the $n \times n$ matrix with a one on the $i$th
diagonal position and zeros otherwise. Hence, the dual
semidefinite problem has the form
\[
\max \sum\limits_{i=1}^n y_i \quad \mbox{s.t.} \quad \diag(y_{1 \;
: \; n})+y_{n+1}(ee^T) \preceq \dfrac{1}{4} L, \quad y \in
\R^{n+1}.
\]
The constraints $\diag(X) = e, \; X \succeq 0$ imply PBQ with
finite upper bounds $\lambda_{\max}(X) \le \ox = n$. Corollary
\ref{coro6.1} yields

\begin{coro} \label{coro9.1}
Let $\ty \in \R^{n+1}$, and assume that the matrix
\[
D = \dfrac{1}{4}L-\Diag(\ty_{1 \; : \; n})-\ty_{n+1}(ee^T)
\]
has at most $l$ negative eigenvalues, and let $\ud \le
\lambda_{\min}(D)$. Then
\[
\hf_p  \ge \sum\limits_{i=1}^n \ty_i + l \cdot n \cdot \ud^ -   =:
\underline{f}_p.
\]
\end{coro}

In Table \ref{tab1} some numerical results for problems given by
Gruber and Rendl \cite{GruRe02} are displayed. The number of nodes
is denoted by $n$. For this suite of ill-posed problems with up to
$601$ constraints and $180 300$ variables the semidefinite
programming solver SDPT3, version 3.02 \cite{TuToTo03} has
computed approximations of the dual optimal value $\tilde{f}_d$,
which are close to the approximate primal optimal value
$\tilde{f}_p$; see the column $\mu(\tilde{f}_p,\tilde{f}_d)$.
Here, the relative accuracy of two real numbers $a$ and $b$ is
measured by the quantity
\[
\mu(a,b):= {\displaystyle \frac{a-b}{\max\{1.0, (|a|+|b|)/2\}}}.
\]
The negative signs in this column show that weak duality is
violated for the computed approximations in four cases. SDPT3 gave
$tc=0$ (normal termination) for the first five ill-posed examples.
Only in the last case $n=600$ the warning $tc=-5$ (that means {\it
: Progress too slow}) was returned. We have computed the lower
bound $\underline{f}_p$ by using Corollary \ref{coro9.1}. The
small quantities $\mu(\tilde{f}_d,\underline{f}_p)$ show that the
overestimation of the rigorous lower bound $\underline{f}_p$ can
be neglected.  In Table \ref{tab1} the times for computing the
approximations with SDPT3, and for computing $\underline{f}_p$ by
using Corollary \ref{coro9.1} are denoted by $t$ and $t_{low}$,
respectively. It follows that the additional time $t_{low}$ for
computing the rigorous bound $\underline{f}_p$ is small compared
to the time $t$ needed for the approximations. This is of the same
tenor as the quotation of Turing at first.

\begin{table}
  \begin{center}
    \begin{normalsize}
      \begin{tabular}{cccccccc}
    \hline
    $n $  & $\underline{f}_p$ & $\mu(\tilde{f}_p,\tilde{f}_d)$ &
    $\mu(\tilde{f}_d,\underline{f}_p)$ & $t $ & $t_{low}$ &
   $tc$
     \\
    \hline

100 & -3.58065e+003 & -7.117e-008 & 3.843e-011 &   4.2 &   0.5 &   0  \\
200 & -1.04285e+004 & -7.018e-008 & 9.621e-010 &   7.9 &   0.2 &    0  \\
300 & -1.90966e+004 & -2.573e-008 & 8.918e-009 &  21.1 &   0.9 &    0  \\
400 & -3.01393e+004 & -1.633e-008 & 3.008e-008 &  39.0 &   2.0 &   0  \\
500 & -4.22850e+004 & 1.431e-008 & 2.584e-008 &  67.5 &   3.7 &   0  \\
600 & -5.57876e+004 & 5.418e-009 & 1.829e-008 & 124.7 &   6.0 &   -5  \\

    \hline
      \end{tabular}
    \end{normalsize}
  \end{center}
  \caption{Results for Graph Partitioning}
  \label{tab1}
\end{table}

\section{Numerical Results}
In this section we describe briefly our numerical experience.
Lurupa \cite{Keil05} is a C++ implementation of the presented
rigorous bounds for the special case of linear programming. In the
following we give a short summary of numerical results for the
NETLIB suite of linear programming problems \cite{Netlib}. For
details refer to \cite{KeilJa04}. The NETLIB LP-suite is a
well-known collection of difficult to solve problems with up to
15695 variables and 16675 constraints. They originate from various
applications, for example forestry, flap settings on aircraft, and
staff scheduling. We chose the set of problems for which
Ord\'o\~nez and Freund \cite{OrdFre03} have computed condition
numbers. There it is stated that 71\% of the problems have an
infinite condition number. As Fourer and Gay \cite{FouGa94}
observed, preprocessing can change the state of an LP from
feasible to infeasible and vice versa, and therefore preprocessing
was not applied.

Roughly speaking, a finite lower bound (upper bound) of the
optimal value can be computed iff the distance to dual
infeasibility (primal infeasibility) is greater than zero. For 76
problems a finite lower bound could be computed with a median
accuracy of $\mbox{med}(\mu(\tilde{f},\overline{f}_p)) = 2.2 \cdot
10^{-8}$ ($\tilde{f}$ is the approximate optimal value) and a
median time ratio of $\mbox{med}(t_{low}/t) = 0.5$. For 35
problems Lurupa has computed a finite upper bound with
$\mbox{med}(\mu(\tilde{f},\overline{f}_d)) = 8.0 \cdot 10^{-9}$
and $\mbox{med}(t_{up}/t) = 5.3$. For 32 well-posed problems
finite rigorous lower and upper bounds could be computed with
$\mbox{med}(\mu(\overline{f}_p,\underline{f}_d)) = 5.6 \cdot
10^{-8}$. Only for two problems with finite condition number
(SCSD8 and SCTAP1) an upper bound of the optimal value could not
be computed. Taking into account the approxmimate solver's default
stopping tolerance of $10^{-9}$, the guaranteed accuracy computed
with Lurupa for the NETLIB LP suite is reasonable. The upper bound
is more expensive, since linear systems have to be solved
rigorously, and sometimes perturbed problems have .

For the SDPLIB benchmark problems of Borchers \cite{Bor99} we have
comuted with VSDP \cite{Ja06a}, a MATLAB software package for
verified semidefinite programming, rigorous bounds . For details
see \cite{Ja06a} and \cite{JaChaKeil05}. Freund,
Ord{\'{o}}{\~{n}}ez and Toh \cite{FreOrdToh06} have solved $85$
problems with SDPT3 out of the $92$ problems of the SDPLIB. They
have omitted the four infeasible problems and three very large
problems where SDPT3 produced out of memory. In their paper
interior-point iteration counts with respect to different measures
for semidefinite programming problems are investigated, and it is
pointed out that $32$ are ill-posed. VSDP could compute (by using
SDPT3 as approximate solver) for all $85$ problems a rigorous
lower bound of the optimal value and could verify the existence of
strictly dual feasible solutions, which proves that all problems
have a zero duality gap. A finite rigorous upper bound could be
computed for all well-posed problems with one exception; this is
{\tt hinf2}. For all $32$ ill-posed problems VSDP has computed the
upper bound $\overline{f}_d = +\infty$, which reflects exactly
that the distance to the next primal infeasible problem is zero as
well as the infinite condition number.

For the $85$ test problems (not counting the $4$ infeasible ones)
SDPT3 (with default values) has computed good approximations and
gave $7$ warnings, and $2$ warnings are given for well-posed
problems. Hence, no warnings are given for $27$ ill-posed problems
with zero distance to primal infeasibility. In other words, there
is no correlation between warnings and the difficulty of the
problem. At least for this test set our rigorous bounds reflect
the difficulty of the problems much better, and they provide
safety, especially in the case where algorithms subsequently call
other algorithms, as is done for example in branch-and-bound
methods.

Some major characteristics of the numerical results of VSDP for
the well-posed SDPLIB-problems are as follows: The median of the
time ratio for computing the rigorous lower (upper) bound and the
approximation is $0.085, (1.99)$, respectively. The median of the
guaranteed accuracy for the problems with finite condition number
is $7.01 \cdot 10^{-7}$. We have used the median here because
there are some outliers.

One of the largest problems which could be rigorously solved by
VSDP is {\tt thetaG51} where the number of constraints is $m =
6910$, and the dimension of the primal symmetric matrix {\tt X} is
$ s = 1001$ (implying $501 501$ variables). For this problem SDPT3
gave the message out of memory, and we used SDPA \cite{SDPA-M} as
approximate solver. The rigorous lower and upper bounds computed
by VSDP are $\underline{f}_p=-3.4900 \cdot 10^2$, $\overline{f}_d=
-3.4406 \cdot 10^2$, respectively. This is an outlier because the
guaranteed relative accuracy is only $0.014$, which may be
sufficient in several applications, but is insufficient from a
numerical point of view. However, existence of optimal solutions
and strong duality is proved. The times in seconds for computing
the approximations, the lower and the upper bound of the optimal
value are {\tt t}$= 3687.95$, $t_{low}= 45.17$, and $t_{up}=
6592.52$, respectively.

For further numerical results and applications concerning
ill-posed problems and the problem of computing the ground state
energy of atomic and molecular systems by using a variational
approach (see for example Fukuda et al. \cite{FuBrNaOvPeYaZh} and
Nakata et al. \cite{NaFuNaFu01}) refer to VSDP \cite{Ja06a}.

To our knowledge no other software packages compute rigorous
results for semidefinite programs. There are several packages that
compute verified results for optimization problems where the
objective and the constraints are defined by smooth algebraic
expressions. Elaborate comparisons  with some of these packages in
the case of linear programming problems can be found in the
forthcoming paper of Keil \cite{Keil2007}.

\section{Conclusions}
The computation of rigorous error bounds for conic optimization
problems can be viewed as a carefully postprocessing tool that
uses only approximate solutions computed by any conic solver. The
bounds are developed in the framework of functional analysis.
Error bounds for special conic problems can be derived easily.

Several numerical results demonstrate that  rigorous error bounds
can be reasonably easily computed even for problems of large size
and for ill-conditioned problems, in most cases with a range which
is not much larger than necessary.

\bibliography{H:/bib/extern,H:/bib/ti3}
\end{document}